\documentclass[a4paper, reqno, 12pt]{amsart}

\usepackage[usenames,dvipsnames]{color}
\usepackage{amsthm,amsfonts,amssymb,amsmath,amsxtra}
\usepackage[all]{xy}
\SelectTips{cm}{}
\usepackage{xr-hyper}
\usepackage[colorlinks=
   citecolor=Black,
   linkcolor=Red,
   urlcolor=Blue]{hyperref}
\usepackage{verbatim}

\usepackage[margin=1.25in]{geometry}
\usepackage{mathrsfs}

\RequirePackage{xspace}
\RequirePackage{etoolbox}
\RequirePackage{varwidth}
\RequirePackage{enumitem}
\RequirePackage{tensor}
\RequirePackage{mathtools}
\RequirePackage{longtable}
\RequirePackage{multirow}

\def\ge{\geqslant}
\def\le{\leqslant}

\def\G{\Gamma}

\def\e{\epsilon}

\def\s{\sigma}

\def\k{\kappa}
\def\l{\lambda}

\def\i{^{-1}}

\def\<{\langle}
\def\>{\rangle}

\newcommand{\BF}{\ensuremath{\mathbb {F}}\xspace}
\newcommand{{\BG}}{\ensuremath{\mathbb {G}}\xspace}

\newcommand{{\BK}}{\ensuremath{\mathbb {K}}\xspace}

\newcommand{\BN}{\ensuremath{\mathbb {N}}\xspace}

\newcommand{\BQ}{\ensuremath{\mathbb {Q}}\xspace}

\newcommand{\BZ}{\ensuremath{\mathbb {Z}}\xspace}

\newcommand{\CO}{\ensuremath{\mathcal {O}}\xspace}

\DeclareMathOperator{\Adm}{Adm}

\def\tW{\tilde W}

\def\kk{\mathbf k}

\newtheorem{theorem}{Theorem}

\theoremstyle{definition}

\numberwithin{equation}{section}
\numberwithin{theorem}{section}

\setitemize[0]{leftmargin=*,itemsep=\the\smallskipamount}
\setenumerate[0]{leftmargin=*,itemsep=\the\smallskipamount}

\renewcommand{\to}{%
   \ifbool{@display}{\longrightarrow}{\rightarrow}%
   }
\let\shortmapsto\mapsto
\renewcommand{\mapsto}{%
   \ifbool{@display}{\longmapsto}{\shortmapsto}%
   }
\newlength{\olen}
\newlength{\ulen}
\newlength{\xlen}
\newcommand{\xra}[2][]{%
   \ifbool{@display}%
      {\settowidth{\olen}{$\overset{#2}{\longrightarrow}$}%
       \settowidth{\ulen}{$\underset{#1}{\longrightarrow}$}%
       \settowidth{\xlen}{$\xrightarrow[#1]{#2}$}%
       \ifdimgreater{\olen}{\xlen}%
          {\underset{#1}{\overset{#2}{\longrightarrow}}}%
          {\ifdimgreater{\ulen}{\xlen}%
             {\underset{#1}{\overset{#2}{\longrightarrow}}}
             {\xrightarrow[#1]{#2}}}}%
      {\xrightarrow[#1]{#2}}
   }
\makeatother
\newcommand{\xyra}[2][]{%
   \settowidth{\xlen}{$\xrightarrow[#1]{#2}$}%
   \ifbool{@display}%
      {\settowidth{\olen}{$\overset{#2}{\longrightarrow}$}%
       \settowidth{\ulen}{$\underset{#1}{\longrightarrow}$}%
       \ifdimgreater{\olen}{\xlen}%
          {\mathrel{\xymatrix@M=.12ex@C=3.2ex{\ar[r]^-{#2}_-{#1} &}}}%
          {\ifdimgreater{\ulen}{\xlen}%
             {\mathrel{\xymatrix@M=.12ex@C=3.2ex{\ar[r]^-{#2}_-{#1} &}}}
             {\mathrel{\xymatrix@M=.12ex@C=\the\xlen{\ar[r]^-{#2}_-{#1} &}}}}}%
      {\mathrel{\xymatrix@M=.12ex@C=\the\xlen{\ar[r]^-{#2}_-{#1} &}}}%
   }
\makeatletter
\newcommand{\xla}[2][]{%
   \ifbool{@display}%
      {\settowidth{\olen}{$\overset{#2}{\longleftarrow}$}%
       \settowidth{\ulen}{$\underset{#1}{\longleftarrow}$}%
       \settowidth{\xlen}{$\xleftarrow[#1]{#2}$}%
       \ifdimgreater{\olen}{\xlen}%
          {\underset{#1}{\overset{#2}{\longleftarrow}}}%
          {\ifdimgreater{\ulen}{\xlen}%
             {\underset{#1}{\overset{#2}{\longleftarrow}}}
             {\xleftarrow[#1]{#2}}}}%
      {\xleftarrow[#1]{#2}}
   }
\newcommand{\isoarrow}{%
   \ifbool{@display}{\overset{\sim}{\longrightarrow}}{\xrightarrow\sim}%
   }
   
\begin{document}

\title[]{Some results on affine Deligne-Lusztig varieties}
\author[X. He]{Xuhua He}
\address{Department of Mathematics, University of Maryland, College Park, MD 20742}
\email{xuhuahe@math.umd.edu}

\thanks{X. H. was partially supported by NSF DMS-1463852}


\begin{abstract}
The study of affine Deligne-Lusztig varieties  originally arose from arithmetic geometry, but many problems on affine Deligne-Lusztig varieties are purely Lie-theoretic in nature. This survey deals with recent progress on several important problems on affine Deligne-Lusztig varieties. The emphasis is on the Lie-theoretic aspect, while some connections and applications to arithmetic geometry will also be mentioned. 
\end{abstract}

\maketitle

\section{Introduction}

\subsection{Bruhat decomposition and conjugacy classes}\label{1.1}
Let $\BG$ be a connected reductive group over a field $\kk$ and $G=\BG(\kk)$. In this subsection, we assume that $\kk$ is algebraically closed. Let $B$ be a Borel subgroup of $G$ and $W$ be the finite Weyl group of $G$. The Bruhat decomposition $G=\sqcup_{w \in W} B w B$ plays a fundamental role in Lie theory. This is explained by Lusztig \cite{Lu-Bruhat} in the memorial conference of Bruhat: 

{\it ``By allowing one to reduce many questions about $G$ to questions about the Weyl group $W$, Bruhat decomposition is indispensable for the understanding of both the structure and representations of $G$.''}

Below we mention two examples of the interaction between the Bruhat decomposition and the (ordinary and twisted) conjugation action of $G$. 

\begin{enumerate}
\item Assume that $\kk=\overline \BF_q$ and $\s$ is the Frobenius of $\kk$ over $\BF_q$. We assume that $\BG$ is defined over $\BF_q$ and we denote by $\s$ the corresponding Frobenius morphism on $G$. The (classical) Deligne-Lusztig varieties was introduced by Deligne and Lusztig in their seminal work \cite{DL}. For any element $w \in W$, the corresponding Deligne-Lusztig variety $X_w$ is a subvariety of the flag variety $G/B$ defined by $$X_w=\{g B \in G/B; g \i \s(g) \in B w B\}.$$

By Lang's theorem, the variety $X_w$ is always nonempty. It is a locally closed, smooth variety of dimension $\ell(w)$. The finite reductive group $\BG(\BF_q)$ acts naturally on $X_w$ and on the cohomology of $X_w$. The Deligne-Lusztig variety $X_w$ plays a crucial role in the representation theory of finite reductive groups, see \cite{DL} and \cite{L-orange}. The structure of $X_w$ has also found important applications in number theory, e.g., in the work of Rapoport, Terstiege and Zhang \cite{RTZ}, and in the work of Li and Zhu \cite{LZ} on the proof of special cases of the ``arithmetic fundamental lemma'' of Zhang \cite{Zhang}.

\item Let $\kk$ be any algebraically closed field. In a series of papers \cite{L1-3}, Lusztig discovered a deep relation between the unipotent conjugacy classes of $G$ and the conjugacy classes of $W$, via the study of the intersection of the unipotent conjugacy classes with the Bruhat cells of $G$. 
\end{enumerate} 

\subsection{Affine Deligne-Lusztig varieties}
The main objects of this survey are affine Deligne-Lusztig varieties, analogous of classical Deligne-Lusztig varieties for loop groups. 

Unless otherwise stated, in the rest of this survey we assume that $\kk=\overline \BF_q((\e))$.  Let $\s$ be the Frobenius morphism of $\kk$ over $\BF_q((\e))$. We assume that $\BG$ is defined over $\BF_q((\e))$ and we denote by $\s$ the corresponding Frobenius morphism on the loop group $G=\BG(\kk)$. We choose a $\s$-stable Iwahori subgroup $I$ of $G$. If $G$ is unramified, then we also choose a $\s$-stable hyperspecial parahoric subgroup $K \supset I$. The affine flag variety $Fl=G/I$ and the affine Grassmannian $Gr=G/K$ (if $G$ is unramified) have natural scheme structures.\footnote{One may replace $\overline \BF_q((\e))$ by the fraction field of the Witt ring. In that case, the affine Grassmannian $Gr$ and the affine flag variety $Fl$ have the structure of perfect schemes, thanks to the recent breakthrough of Zhu \cite{Z2}, and of Bhatt and Scholze \cite{BS}. Many of the results we discuss in this survey hold for the fraction field of the Witt ring as well.}

Let $S$ be a maximal $\kk$-split torus of $G$ defined over $\BF_q((\e))$ and let $T$ be its centralizer, a maximal torus of $G$. The Iwahori-Weyl group associated to $S$ is $$\tW=N(\kk)/T(\kk)_1,$$ where $N$ is the normalizer of $S$ in $G$ and $T(\kk)_1$ is the maximal open compact subgroup of $T(\kk)$. The group $\tW$ is also a split extension of the relative (finite) Weyl group $W_0$ by the normal subgroup $X_*(T)_{\G_0}$, where $X_*(T)$ is the coweight lattice of $T$ and $\G_0$ is the Galois group of $\overline{\kk}$ over $\kk$  (cf. \cite[Appendix]{PR}). The group $\tW$ has a natural quasi-Coxeter structure. We denote by $\ell$ and $\le$ the length function and the Bruhat order on $\tW$. We have the following generalization of the Bruhat decomposition $$G=\sqcup_{w \in \tW} I w I,$$ due to Iwahori and Matsumoto \cite{IM} in the split case, and to Bruhat and Tits \cite{BT} in the general case. If $G$ is unramified, then we also have $$G=\sqcup_{\l \text{ is a dominant coweight}} K \e^\l K.$$ 

Affine Deligne-Lusztig varieties were introduced by Rapoport in \cite{Ra}. Compared to the classical Deligne-Lusztig varieties, we need two parameters here: an element $w$ in the Iwahori-Weyl group $W$ and an element $b$ in the loop group $G$. The corresponding affine Deligne-Lusztig variety (in the affine flag variety) is defined as $$X_w(b)=\{g I \in G/I; g \i b \s(g) \in I w I\} \subset Fl.$$
If $G$ is unramified, one may use a dominant coweight $\l$ instead of an element in $\tW$ and define the affine Deligne-Lusztig variety (in the affine Grassmannian) by $$X_\l(b)=\{g K \in G/K; g \i b \s(g) \in K \e^\l K\} \subset Gr.$$

Affine Deligne-Lusztig varieties are schemes locally of finite type over $\overline \BF_q$. Also the varieties are isomorphic if the element $b$ is replaced by another element $b'$ in the same $\s$-conjugacy class. 

A major difference between affine Deligne-Lusztig varieties and  classical Deligne-Lusztig varieties is that affine Deligne-Lusztig varieties have the second parameter: the element $b$, or the $\s$-conjugacy class $[b]$ in the loop group $G$; while in the classical case considered in \S \ref{1.1}, by Lang's theorem there is only one $\s$-conjugacy class in $\BG(\overline \BF_q)$ and thus adding a parameter $b \in \BG(\overline \BF_q)$ does not give any new variety. 

The second parameter $[b]$ in the affine Deligne-Lusztig varieties makes them rather challenging to study, both from the Lie-theoretic point of view, and from the arithmetic-geometric point of view. Below we list some major problems on the affine Deligne-Lusztig varieties:

\begin{itemize}
\item When is an affine Deligne-Lusztig variety nonempty?

\item If it is nonempty, what is its dimension?

\item What are the connected components? 

\item Is there a simple geometric structure for certain affine Deligne-Lusztig varieties?

\end{itemize}

We may also consider the affine Deligne-Lusztig varieties associated to arbitrary parahoric subgroups, besides hyperspecial subgroups and Iwahori subgroups. This will be discussed in \S \ref{par}. 

\subsection{A short overview of $X(\mu, b)$}\label{union}
The above questions may also be asked for a certain union $X(\mu, b)$ of affine Deligne-Lusztig varieties in the affine flag variety. 

Let $\mu$ be a dominant coweight of $G$ with respect to a given Borel subgroup of $G$ over $\kk$ (in applications to number theory, $\mu$ usually comes from a Shimura datum). The admissible set $\Adm(\mu)$ was introduced by Kottwitz and Rapoport in \cite{KR}. It is defined by $$\Adm(\mu)=\{w \in \tW; w \le t^{x(\mu)} \text{ for some } x \in W_0\}.$$ 

We may explain it in a more Lie-theoretic language. Let $Gr_{\mathcal G}$ be the deformation from the affine Grassmannian to the affine flag variety \cite{Gai}. The coherence conjecture of Pappas and Rapoport \cite{PR} implies that the special fiber of the global Schubert variety $\overline{Gr}_{\mathcal G, \mu}$ associated to the coweight $\mu$ (cf. \cite[Definition 3.1]{Z1}) is $\cup_{w \in \Adm(\mu)} I w I/I$. This conjecture was proved by Zhu in \cite{Z1}. Now we set $$X(\mu, b)=\cup_{w \in \Adm(\mu)} X_w(b) \subset Fl.$$ This is a closed subscheme of $Fl$ and serves as the group-theoretic model for the Newton stratum corresponding to $[b]$ in the special fiber of a Shimura variety giving rise to the datum $(G, \mu)$. 

It is also worth mentioning that, although the admissible set $\Adm(\mu)$ has a rather simple definition, it is a very complicated combinatorial object. We refer to the work of Haines and Ng\^{o} \cite{HN}, and the recent joint work of the author with Haines \cite{HH} for some properties of $\Adm(\mu)$. 

\subsection{Current status} Affine Deligne-Lusztig varieties in the affine Grassmannian are relatively more accessible than the ones in the affine flag variety, mainly due to the following two reasons: 
\begin{itemize} 
\item The set of dominant coweights is easier to understand than the Iwahori-Weyl group; 

\item For $X_\l(b)$ the group $G$ is unramified while for $X_w(b)$, we need to deal with ramified, or even non quasi-split reductive groups. 
\end{itemize}

For an unramified group $G$, we also have the fibration $\cup_{w \in W_0 t^\l W_0} X_w(b) \to X_{\l}(b)$, with fibers isomorphic to the flag variety of $\BG(\overline \BF_q)$. Thus much information on $X_\l(b)$ can be deduced from $X_w(b)$. 

Nevertheless, the study of the affine Deligne-Lusztig varieties in affine Grassmannian is a very challenging task and has attracted the attention of experts in arithmetic geometry in the past two decades. It is a major achievement in arithmetic geometry to obtain a fairly good understanding on these varieties. 

As to the affine Deligne-Lusztig varieties in the affine flag varieties, the situation is even more intriguing. We have made significant progress in the past 10 years in this direction, yet many aspects of $X_w(b)$ remain rather mysterious. I hope that by combining various Lie-theoretic methods together with arithmetic-geometric methods, our knowledge on affine Deligne-Lusztig varieties will be considerably advanced. 

In the rest of the survey, we will report on some recent progress on the affine Deligne-Lusztig varieties. 

\subsection*{Acknowledgement} We thank Ulrich G\"ortz, Urs Hartl, George Lusztig, Michael Rapoport, Sian Nie and Rong Zhou for useful comments. 

\section{Some relation with affine Hecke algebras}

\subsection{The set $B(G)$ and Kottwitz's classification}
Let $B(G)$ be the set of $\s$-conjugacy classes of $G$.  Kottwitz \cite{K1} and \cite{K2} gave a classification of the set $B(G)$, generalizing the Dieudonn\'e-Manin classification of isocrystals by their Newton polygons. Any $\s$-conjugacy class $[b]$ is determined by two invariants: 
\begin{itemize}
\item The element $\k([b]) \in \pi_1(G)_{\G}$, where $\G$ is the Galois group of $\overline \kk$ over $\BF_q((\e))$; 

\item The Newton point $\nu_b$ in the dominant chamber of $X_*(T)_{\G_0} \otimes \BQ$. 
\end{itemize}

A different point of view, which is quite useful in this survey, is the relation between the set $B(G)$ with the set $B(\tW, \s)$ of $\s$-conjugacy classes of $\tW$. Recall that $\tW=N(\kk)/T(\kk)_1$. The natural embedding $N(\kk) \to G$ induces a natural map $\Psi: B(\tW, \s) \to B(G)$. By \cite{GHKR2} and \cite{He14}, the map $\Psi$ is surjective. The map $\Psi$ is not injective. However, there exists an important family $B(\tW, \s)_{str}$ of straight $\s$-conjugacy classes of $\tW$. By definition, a $\s$-conjugacy class $\CO$ of $\tW$ is {\it straight} if it contains an element $w \in \CO$ such that $\ell(w \s(w) \cdots \s^{n-1}(w))=n \ell(w)$ for all $n \in \BN$. 
The following result is discovered in \cite[Theorem 3.7]{He14}. 
\begin{theorem}
The map  $\Psi: B(\tW, \s) \to B(G)$ induces a bijection  $$B(\tW, \s)_{str} \longleftrightarrow B(G).$$
\end{theorem}

This result gives the parametrization of the $\s$-conjugacy classes of $G$ in terms of the set of straight $\s$-conjugacy classes of its Iwahori-Weyl group $\tW$. In particular, the two parameters occurring in the definition of the affine Deligne-Lusztig variety $X_w(b)$ are all from $\tW$. 


\smallskip

Note that the affine Deligne-Lusztig variety $X_w(b)$ is closely related to the intersection $I w I \cap [b]$. This intersection is very complicated in general. However, it is discovered in \cite{He14} that for certain elements $w \in \tW$, the intersection $I w I \cap [b]$ equals $I w I$. More precisely, we denote by $\tW_{\s-\min}$ the set of elements in $\tW$ that are of minimal length in their $\s$-conjugacy classes. Then
$$\text{ For $w \in \tW_{\s-\min}$, $I w I  \subset [b]$ if $[b]=\Psi(w)$.}$$ This serves as the starting point of the reduction method for affine Deligne-Lusztig varieties $X_w(b)$ for arbitrary $w$. 

\subsection{``Dimension=Degree'' theorem}
Deligne and Lusztig introduced in \cite{DL} a reduction method to study the classical Deligne-Lusztig varieties. Their method works for the affine Deligne-Lusztig varieties as well. Some  combinatorial properties of affine Weyl groups established in joint work with Nie \cite{HN14} allow us to reduce the study of $X_w(b)$ for any $w$, via the reduction method \`a la Deligne and Lusztig, to the study of $X_w(b)$ for $w \in \tW_{\s-\min}$. 

The explicit reduction procedure, however, is rather difficult to keep track of. In \cite{He14}, we discovered that the same reduction procedure appears in a totally different context as follows. 

Let $H$ be the affine Hecke algebra (over $\BZ[v^{\pm 1}]$) associated to $\tW$. Let $[\tilde H, \tilde H]_{\s}$ be the $\s$-twisted commutator, i.e. the $\BZ[v^{\pm 1}]$-submodule of $H$ generated by $[h, h']_\s=h h'-h'\s(h)$. By \cite{HN14}, the $\s$-twisted cocenter $\overline{H}=H/[H, H]_\s$ has a standard basis given by $\{T_{\CO}\}$, where $\CO$ runs over all the $\s$-conjugacy classes of $\tW$. Thus for any $w \in \tW$, we have $$T_w \equiv \sum_{\CO} f_{w, \CO} T_{\CO} \mod [H, H]_\s.$$ The coefficients $f_{w, \CO} \in \BN[v-v \i]$, which we call the class polynomials  (over $v-v \i$). We have the following ``dimension=degree'' theorem established in \cite{He14}. 

\begin{theorem}\label{deg=dim}
Let $b \in G$ and $w \in \tW$. Then $$\dim(X_w(b))=\max_{\CO; \Psi(\CO)=[b]} \frac{1}{2}\bigl(\ell(w)+\ell(\CO)+\deg(f_{w, \CO}) \bigr)-\<\nu_b, 2 \rho\>.$$
\end{theorem}

Here $\ell(\CO)$ is the length of any minimal length element in $\CO$ and $\rho$ is the half sum of positive roots in $G$. Here we use the convention that the dimension of an empty variety and the degree of a zero polynomial are both $-\infty$. Thus the above theorem reduces the nonemptiness question and the dimension formula of $X_w(b)$ to some questions on the class polynomials $f_{w, \CO}$ for $\Psi(\CO)=[b]$. 

The explicit computation of the class polynomials is very difficult at present. Note that there is a close relation between the cocenter and representations of affine Hecke algebras \cite{CH}. One may hope that some progress in the representation theory of affine Hecke algebras would also advance our knowledge on affine Deligne-Lusztig varieties. At present, we combine the ``dimension=degree'' theorem together with some Lie-theoretic techniques, and the results on $X_\l(b)$ in the affine Grassmannian established previously by arithmetic-geometric method, to obtain some explicit answers to certain questions on $X_w(b)$ and on $X(\mu, b)$. 

\section{Nonemptiness pattern}
\subsection{Mazur's inequality} In this subsection, we discuss the  non-emptiness patterns of affine Deligne-Lusztig varieties. Here Mazur's inequality plays a crucial role. 

In \cite{Ma}, Mazur proved that the Hodge slope of any $F$-crystal is always larger than or equal to the Newton slope of associated isocrystal. The converse was obtained by Kottwitz and Rapoport in \cite{KR03}. Here we regard the Newton slope and Hodge slope as elements in $\BQ^n_+=\{a_1, \cdots, a_n; a_1 \ge \cdots \ge a_n\}$ and the partial order in $\BQ^n_+$ is the dominance order, i.e. $(a_1, \cdots, a_n) \preceq (b_1, \cdots, b_n)$ if and only if $a_1 \le b_1, a_1+a_2 \le b_1+b_2, \cdots, a_1+\cdots+a_{n-1} \le b_1+\cdots+b_{n-1}, a_1+\cdots+a_n=b_1+\cdots+b_n$. 

Note that $\BQ^n_+$ is the set of rational dominant coweights for $GL_n$. The dominant order can be defined for the set of rational dominant coweights for any reductive group. This is what we use to describe the nonemptiness pattern of some affine Deligne-Lusztig varieties. 

\subsection{In the affine Grassmannian} For $X_\l(b)$ in the affine Grassmannian, we have a complete answer to the nonemptiness question. 

\begin{theorem}\label{3.1}
Let $\l$ be a dominant coweight and $b \in G$. Then $X_\l(b) \neq \emptyset$ if and only if $\k([b])=\k(\l)$ and $\nu_b \preceq \l$. 
\end{theorem}
The ``only if'' part was proved by Rapoport and Richartz in \cite{RR}, and by Kottwitz in \cite{K3}. The ``if'' part was proved by Gashi \cite{Ga}. The result also holds if the hyperspecial subgroup of an unramified group is replaced by a maximal special parahoric subgroup of an arbitrary reductive group. This was obtained in \cite{He14} using the ``dimension=degree'' Theorem \ref{deg=dim}. 

\subsection{In the affine flag}\label{0.7.2} Now we consider the variety $X_w(b)$ in the affine flag variety. 

(i) We first discuss the case where $[b]$ is basic, i.e., the corresponding Newton point $\nu_b$ is central in $G$ (and thus Mazur's inequality is automatically satisfied). 

\begin{theorem}
Let $G$ be a quasi-split group. Let $[b] \in B(G)$ be basic and $w \in \tW$. Then $X_w(b) \neq \emptyset$ if and only if there is no ``Levi obstruction''. 
\end{theorem}

The ``Levi obstruction'' is defined in terms of the $P$-alcove elements, introduced by G\"ortz, Haines, Kottwitz, and Reuman in \cite{GHKR2}. The explicit definition is technical and we omit it here. This result was conjectured by G\"ortz, Haines, Kottwitz, and Reuman in \cite{GHKR2} for split groups and was established in joint work with G\"ortz and Nie \cite{GHN} for any quasi-split group. Note that the ``quasi-split'' assumption here is not essential as one may relate $X_w(b)$ for any reductive group $G$ to another affine Deligne-Lusztig variety for the quasi-split inner form of $G$. We refer to \cite[Theorem 2.27]{CDM} for the explicit statement in the general setting. 

\smallskip

(ii) For any nonbasic $\s$-conjugacy class $[b]$, one may ask for analogues of ``Mazur's inequality'' and/or the ``Levi obstruction'' in order to describe the nonemptiness pattern of $X_w(b)$. This is one of the major open problems in this area. We refer to \cite[Remark 12.1.3]{GHKR2} for some discussion in this direction. As a first step, one may consider the conjecture of G\"ortz-Haines-Kottwitz-Reumann \cite[Conjecture 9.5.1 (b)]{GHKR2} on the asymptotic behavior of $X_w(b)$ for nonbasic $[b]$. Some affirmative answer to this conjecture was given in \cite[Theorem 2.28]{CDM} and \cite{MST} in the case where $[b]=[\e^\l]$ for some dominant coweight $\l$. 


\subsection{Kottwitz-Rapoport conjecture} To describe the nonemptiness pattern on the union $X(\mu, b)$  of affine Deligne-Lusztig varieties in the affine flag variety, we recall the definition of neutrally acceptable $\s$-conjugacy classes introduced by Kottwitz in \cite{K2}, $$B(G, \mu)=\{[b] \in B(G); \k([b])=\k(\mu), \nu_b \le \mu^\diamond\},$$ where $\mu^\diamond$ is the Galois average of $\mu$. 

By Theorem \ref{3.1}, $X_\mu(b) \neq \emptyset$ if and only if $[b] \in B(G, \mu)$. We have a similar result for the union $X(\mu, b)$ of affine Deligne-Lusztig varieties in the affine flag variety.

\begin{theorem}\label{KR-conj}
Let $[b] \in B(G)$. Then $X(\mu, b) \neq \emptyset$ if and only if $[b] \in B(G, \mu)$. 
\end{theorem}

This result was conjectured by Kottwitz and Rapoport in \cite{KR03} and \cite{Ra}. The ``only if'' part is a group-theoretic version of Mazur's inequality and was proved by Rapoport and Richartz for unramified groups in \cite[Theorem 4.2]{RR}. The ``if'' part is the ``converse to Mazur's inequality'' and was proved by Wintenberger in \cite{Wi} for quasi-split groups. The general case in both directions was established in \cite{He16} by a different approach, via a detailed analysis of the map $\Psi: B(\tilde W) \to B(G)$, of the partial orders on $B(G)$ (an analogy of Grothendieck's conjecture for the loop groups) and of the maximal elements in $B(G, \mu)$ \cite{HNx}.  

As we mentioned in \S\ref{0.7.2}, for a single affine Deligne-Lusztig variety $X_w(b)$, one may reduce the case of a general group to the quasi-split case. However, for the union of affine Deligne-Lusztig varieties, the situation is different. There is no relation between the admissible set $\Adm(\mu)$ (and hence $X(\mu, b)$) for an arbitrary reductive group and its quasi-split inner form. This adds essential difficulties in the study of $X(\mu, b)$ for non quasi-split groups.  

Rad and Hartl in \cite{RH} established the analogue of the Langlands-Rapoport conjecture \cite{LR} for the rational points in the moduli stacks of global $G$-shtukas, for arbitrary connected reductive groups and arbitrary parahoric level structure. They described the rational points as a disjoint union over isogeny classes of global $G$-Shtukas, and then used Theorem \ref{KR-conj} to determine which isogeny classes are nonempty. 

\section{Dimension formula}

\subsection{In the affine Grassmannian} For $X_\l(b)$ in the affine Grassmannian, we have an explicit dimension formula. 

\begin{theorem}\label{4.1}
Let $\l$ be a dominant coweight and $b \in G$. If $X_\l(b) \neq \emptyset$, then $$\dim X_\l(b)=\<\l-\nu_b, \rho\>-\frac{1}{2} \text{def}_G(b),$$ where $\text{def}_G(b)$ is the defect of $b$. 
\end{theorem}

The dimension formula of $X_{\l}(b)$ was conjectured by Rapoport in \cite{Ra}, inspired by Chai's work \cite{Ch}. The current reformulation is due to Kottwitz \cite{K4}. For split groups, the conjectural formula was obtained by G\"ortz, Haines, Kottwitz and Reuman \cite{GHKR1} and Viehmann \cite{Vi06}. The conjectural formula for general quasi-split unramified groups was obtained independently by Zhu \cite{Z2} and Hamacher \cite{Ha1}. 

\subsection{In the affine flag variety} Now we consider $X_w(b)$ in the affine flag variety. 

\begin{theorem}
Let $[b] \in B(G)$ be basic and $w \in \tW$ be an element in the shrunken Weyl chamber (i.e., the lowest two-sided cell of $\tW$). If $X_w(b) \neq \emptyset$, then $$\dim X_w(b)=\frac{1}{2}(\ell(w)+\ell(\eta_\s(w))-\text{def}_G(b)).$$ Here $\eta_\s: \tW \to W_0$ is defined in \cite{GHKR2}.
\end{theorem}

This dimension formula was conjectured by G\"ortz, Haines, Kottwitz, and Reuman in \cite{GHKR2} for split groups and was established for residually split groups in \cite{He14}. The proof in \cite{He14} is based on the ``dimension=degree'' Theorem \ref{deg=dim}, some results on the $\s$-twisted cocenter $\overline{H}$ of affine Hecke algebra $H$, together with the dimension formula of $X_\l(b)$ (which was only known for split groups at that time).  The dimension formula for arbitrary reductive groups (under the same assumption on $b$ and $w$) is obtained by the same argument in \cite{He14}, once the dimension formula of $X_\l(b)$ for quasi-split unramified groups became available, cf. Theorem \ref{4.1}. 

Note that the assumption that $w$ is contained in the lowest two-sided cell is an essential assumption here. A major open problem is to understand the dimension of $X_w(b)$ for $[b]$ basic, when $w$ is in the critical stripes (i.e., outside the lowest two-sided cell). So far, no conjectural dimension formula has been formulated. However, the ``dimension=degree'' Theorem \ref{deg=dim} and the explicit computation in low rank cases \cite{GHKR2} indicate that this problem might be closely related to the theory of Kazhdan-Lusztig cells. I expect that further progress on the affine cellularity of affine Hecke algebras, which is a big open problem in representation theory, might shed new light on the study of $\dim X_w(b)$.

I also would like to point out that affine Deligne-Lusztig varieties in affine Grassmannians are equi-dimensional, while in general affine Deligne-Lusztig varieties in the affine flag varieties are not equi-dimensional. 

\subsection{Certain unions} 

We will see in \S \ref{nice} that for certain pairs $(G, \mu)$, $X(\mu, b)$ admits some simple geometric structure. In these cases, one may write down an explicit dimension formula for $X(\mu, b)$. Outide these case, very little is known for $\dim X(\mu, b)$. 

Here we mention one difficult case: the Siegel modular variety case. Here $G=Sp_{2 g}$ and $\mu$ is the minuscule coweight. It was studied by G\"ortz and Yu in \cite{GY}, in which they showed that for basic $[b]$, $\dim X(\mu, b)=\frac{g^2}{2}$ if $g$ is even and $\frac{g(g-1)}{2} \le \dim X(\mu, b) \le [\frac{g^2}{2}]$ if $g$ is odd. It would be interesting to determine the exact dimension when $g$ is odd. 

\section{Hodge-Newton decomposition}\label{HN} To study the set-theoretic and geometric properties of affine Deligne-Lusztig varieties, a very useful tool is to reduce the study of affine Deligne-Lusztig varieties of a connected reductive group to certain affine Deligne-Lusztig varieties of its Levi subgroups. Such reduction is achieved by the Hodge-Newton decomposition, which originated in Katz's work \cite{Ka} on $F$-crystals with additional structures. In this section, we discuss its variation for affine Deligne-Lusztig varieties in affine Grassmannians, and further development on affine Deligne-Lusztig varieties in affine flag varieties, and on the union of affine Deligne-Lusztig varieties. 

\subsection{In the affine Grassmannian}\label{0.9.1} For affine Deligne-Lusztig varieties in the affine Grassmannian, Kottwitz in \cite{K3} (see also \cite{Vi08}) established the following Hodge-Newton decomposition, which is the group-theoretic generalization of Katz's result. Here the pair $(\l, b)$ is called Hodge-Newton decomposable with respect to a proper Levi subgroup $M$ if $b \in M$ and $\l$ and $b$ have the same image under the Kottwitz's map $\k_M$ for $M$. 

\begin{theorem}
Let $M$ be a Levi subgroup of $G$ and $(\l, b)$ be Hodge-Newton decomposable with respect to $M$. Then the natural map $X^M_\l(b) \to X^G_\l(b)$ is an isomorphism.
\end{theorem}

\subsection{In the affine flag variety} For affine Deligne-Lusztig varieties in affine flag varieties, the situation is more complicated, as the Hodge-Newton decomposability condition on the pairs $(w, b)$ is rather difficult. As pointed out in \cite{GHKR2}, ``It is striking that the notion of $P$-alcove, discovered in the attempt to understand the entire emptiness pattern for the $X_x(b)$ when $b$ is basic, is also precisely the notion needed for our Hodge-Newton decomposition.'' 

The Hodge-Newton decomposition for $X_w(b)$ was established by G\"ortz, Haines, Kottwitz and Reuman in \cite{GHKR2}.

\begin{theorem}
Suppose that $P=M N$ is a semistandard Levi subgroup of $G$ and $w \in \tW$ is a $P$-alcove element in the sense of \cite{GHKR2}. Let $b \in M$. Then the natural map $X^M_w(b) \to X^G_w(b)$ induces a bijection $$J^M_b \backslash X^M_w(b) \cong J^G_b \backslash X^G_w(b).$$
\end{theorem}

\subsection{Certain unions} For $X(\mu, b)$, the Hodge-Newton decomposability condition is still defined on the pair $(\mu, b)$. However, the precise condition is more complicated than in \S \ref{0.9.1} as we consider arbitrary connected reductive groups, not only the unramified ones. We refer to \cite[Definition 2.1]{GHN2} for the precise definition. The following Hodge-Newton decomposition for $X(\mu, b)$ was established in a joint work with G\"ortz and Nie \cite{GHN2}.

\begin{theorem}\label{HN-GHN}
Suppose that $(\mu, b)$ is Hodge-Newton decomposable with respect to some proper Levi subgroup. Then
\[
X(\mu, b) \cong \bigsqcup_{P'=M'N'} X^{M'}(\mu_{P'}, b_{P'}),
\]
where $P'$ runs through a certain finite set of semistandard parabolic subgroups. The subsets in the union are open and closed.
\end{theorem}

We refer to \cite[Theorem 3.16]{GHN2} for the precise statement. Note that an essential new feature is that unlike the Hodge-Newton decomposition of a single affine Deligne-Lusztig variety (e.g. $X_\l(b)$ or $X_w(b)$) where only one Levi subgroup is involved, in the Hodge-Newton decomposition of $X(\mu, b)$ several Levi subgroups are involved. 

Thus, the statement here is more complicated than the Hodge-Newton decomposition of $X_\l(b)$ and $X_w(b)$. But this is consistent with the fact that the Newton strata in the special fiber of Shimura varieties with Iwahori level structure are more complicated than those with hyperspecial level structure. I believe that the Hodge-Newton decomposition here would help us to overcome some of the difficulties occurring in the study of Shimura varieties with Iwahori level structure (as well as arbitrary parahoric level structures). We will see some results in this direction in \S \ref{conn} and in \S \ref{nice}. 

\section{Connected components}\label{conn} In this subsection, we discuss the set of connected components of some closed affine Deligne-Lusztig varieties, e.g. $$X_{\preceq \l}(b):=\cup_{\l' \preceq \l} X_{\l'}(b) \text{ and } X(\mu, b)=\cup_{w \in \Adm(\mu)} X_w(b).$$ The explicit description of the set of connected components has some important applications in number theory, which we will mention later. 

Note that affine Grassmannians and affine flag varieties are not connected in general, and their connected components  are indexed by $\pi_1(G)_{\G_0}$. This gives the first obstruction to the connectedness. The second obstruction comes from the Hodge-Newton decomposition, which we discussed in \S \ref{HN}. One may expect that these are the only obstructions. We have the following results. 

\begin{theorem}
Assume that $G$ is an unramified simple group and that $(\l, b)$ is Hodge-Newton indecomposable. Then $$\pi_0(X_{\preceq \l}(b)) \cong \pi_1(G)_{\G_0}^\s.$$
\end{theorem}

This was first proved by Viehmann for split groups, and then by Chen, Kisin and Viehmann \cite{CKV} for quasi-split unramified  groups and for $\l$ minuscule. The description of $\pi_0(X_{\preceq \l}(b))$ for $G$ quasi-split unramified, and $\l$ non-minuscule, was conjectured in \cite{CKV} and was established by Nie \cite{Nie}.

Note that the minuscule coweight case is especially important for applications in number theory. Kisin \cite{Ki} proved the Langlands-Rapoport conjecture for mod-$p$ points on Shimura varieties of {\it abelian type} with {\it hyperspecial level structure}. Compared to the function field analogous of Langlands-Rapoport conjecture \cite{RH}, there are extra complication coming from algebraic geometry and the explicit description of the connected components of $X(\mu, b)$ in \cite{CKV} is used in an essential way to overcome the complication. 

\begin{theorem}\label{x-conn}
Let $\mu$ be a dominant coweight and $b \in G$. Assume that $[b] \in B(G, \mu)$ and that $(\mu, b)$ is Hodge-Newton indecomposable. Then 

(1) If $[b]$ is basic, then $\pi_0(X(\mu, b)) \cong \pi_1(G)_{\G_0}^\s$. 

(2) If $G$ is split, then $\pi_0(X(\mu, b)) \cong \pi_1(G)$. 
\end{theorem}

Here part (1) was obtained in joint work with Zhou \cite{HZ}. As an application, we verified the Axioms in \cite{HR} for certain PEL type Shimura varieties. In \cite{HZ}, the set of connected components of $X(\mu, b)$ was also studied for nonbasic $b$.  We proved the in a residually split group, the set of connected components is ``controlled" by the set of straight elements, together with the obstruction from the corresponding Levi subgroup. Combined with the work of Zhou \cite{Zh}, we verified in the residually split case, the description of the mod-$p$ isogeny classes on Shimura varieties conjectured by Langlands and Rapoport \cite{LR}. Part (2) is recent work of Chen and Nie \cite{CN}. 

We would like to point out that in the statement, the following two conditions are essential:
\begin{itemize}
\item The $\s$-conjugacy class $[b]$ is neutrally acceptable, i.e. $[b] \in B(G, \mu)$. This condition comes from the Kottwitz-Rapoport conjecture (see Theorem \ref{KR-conj}). 

\item The pair $(\mu, b)$ is Hodge-Newton indecomposable. In the general case, we need to apply the Hodge-Newton decomposition (see Theorem \ref{HN-GHN}). As a consequence, several $\pi_1(M)$ are involved in the description of $\pi_0(X(\mu, b))$ in general. 
\end{itemize}

\section{Arbitrary parahoric level structure}\label{par}

\subsection{Parahoric level versus Iwahori level} Let $K' \supset I$ be a standard parahoric subgroup of $G$ and $W_{K'}$ be the finite Weyl group of $K'$. We define $$X(\mu, b)_{K'}=\{g K' \in G/K'; g \i b \s(g) \in K' \Adm(\mu) K'\}.$$

If $K'=I$, then $X(\mu, b)_{K'}=X(\mu, b)$. If $G$ is unramified, $\mu$ is minuscule and $K'=K$ is a hyperspecial parahoric subgroup, then $X(\mu, b)_{K'}=X_\mu(b)$. As we have mentioned, the varieties $X_\mu(b)$ (resp. $X(\mu, b)$) serve as group-theoretic models for the Newton strata in the special fiber of Shimura varieties with hyperspecial (resp. Iwahori) level structure. The variety $X(\mu, b)_{K'}$ plays the same role in the study of Shimura varieties with arbitrary parahoric level structure. 

The following result relates $X(\mu, b)_{K'}$ for an arbitrary parahoric subgroup $K'$ with $X(\mu, b)$ (for the Iwahori subgroup $I$). 

\begin{theorem}
The projection map $G/I \to G/K'$ induces a surjection $$X(\mu, b) \twoheadrightarrow X(\mu, b)_{K'}.$$
\end{theorem}

This was conjectured by Kottwitz and Rapoport in \cite{KR2} and \cite{Ra} and was proved in \cite{He16}.  This fact allows one to reduce many questions (e.g. nonemptiness pattern, connected components, etc.) of $X(\mu, b)_{K'}$ for arbitrary $K'$ to the same questions for $X(\mu, b)$. In fact, the statements in Theorem \ref{KR-conj} and Theorem \ref{x-conn} hold if $X(\mu, b)$ is replaced by $X(\mu, b)_{K'}$ for an arbitrary parahoric subgroup $K'$. 

\subsection{Lusztig's $G$-stable pieces} I would like to draw attention to some crucial ingredient in the proof, which has important applications in arithmetic geometry. 

Note that $I \Adm(\mu) I \subsetneqq K' \Adm(\mu) K'$ if $I \subsetneqq K'$. In order to show that $X(\mu, b) \to X(\mu, b)_{K'}$ is surjective, one needs to have some decomposition of $K' \Adm(\mu) K'$, finer than the decomposition into $K'$ double cosets. The idea of the sought-after decomposition is essentially due to Lusztig. In \cite{L1}, Lusztig introduced $G$-stable pieces for reductive groups over algebraically closed fields. The closure relation between $G$-stable pieces was determined in \cite{He07} and a more systematic approach using the `` partial conjugation action'' technique was given later in \cite{He072}. The notion and the closure relation of $G$-stable pieces also found application in arithmetic geometry, e.g. in the work of Pink, Wedhorn and Ziegler on algebraic zip data \cite{PWZ}. 

\subsection{Ekedahl-Kottwitz-Oort-Rapoport stratification} In \cite{L99}, Lusztig extended his ideas to the loop groups, see also \cite{He11} and \cite{Vi2}. It was used it to define the Ekedahl-Oort stratification of a general Shimura variety. 

The desired decomposition of $K' \Adm(\mu) K'$ for an arbitrary parahoric subgroup $K'$ was given in \cite{He16} as $$K' \Adm(\mu) K'=\sqcup_{w \in {}^{K'} \tW \cap \Adm(\mu)} K' \cdot_\s I w I,$$ where ${}^{K'} \tW$ is the set of minimal length elements in $W_{K'} \backslash \tW$ and $\cdot_\s$ means the $\s$-conjugation action. This decomposition is used in  joint work with Rapoport \cite{HR} to define the Ekedahl-Kottwitz-Oort-Rapoport stratification of Shimura varieties with arbitrary parahoric level structure. This stratification interpolates between the Kottwitz-Rapoport stratification in the case of the Iwahori level structure and the Ekedahl-Oort stratification \cite{Vi2} in the case of hyperspecial level structure.

\section{Affine Deligne-Lusztig varieties with simple geometric structure}\label{nice}

\subsection{Simple geometric structure for some $X(\mu, b_0)_{K'}$}\label{8.1} The geometric structure of $X(\mu, b_0)_{K'}$ for basic $b_0$ is rather complicated in general. However, in certain cases, $X(\mu, b_0)_{K'}$ admit a simple description. The first nontrivial example is due to Vollaard and Wedhorn in \cite{VW}. They showed that $X_\mu(b_0)$ for an unramified unitary group of signature $(1, n-1)$ and $\mu=(1, 0, \cdots, 0)$ (and for hyperspecial parahoric level structure),  is a union of classical Deligne-Lusztig varieties, and the index set and the closure relations between the strata are encoded in a Bruhat-Tits building. Since then, this question has attracted significant attention. We mention the work of Rapoport, Terstiege and Wilson \cite{RTW} on ramified unitary groups, of Howard and Pappas \cite{HP2}, \cite{HP} on orthogonal groups, of Tiao and Xiao \cite{TX} in the Hilbert-Blumenthal case. In all these works, the parahoric subgroups involved are hyperspecial parahoric subgroups or certain maximal parahoric subgroups. The analogous group-theoretic question for maximal parahoric subgroups was studied in joint work with G\"ortz \cite{GH}. 

Note that these simple descriptions of closed affine Deligne-Lusztig varieties (and the corresponding basic locus of Shimura varieties) have been used, with great success, towards applications in number theory: to compute intersection numbers of special cycles, as in the Kudla-Rapoport program \cite{KR2} or in work \cite{RTZ}, \cite{LZ} towards Zhang's Arithmetic Fundamental Lemma \cite{Zhang}; and to prove the Tate conjecture for certain Shimura varieties \cite{TX2}, \cite{DTX}. 

The work of \cite{VW}, \cite{RTW}, \cite{HP2}, \cite{HP}, \cite{TX} focused on specific Shimura varieties with certain maximal parahoric level structure. The work \cite{GH} studied the analogous group-theoretic question for arbitrary reductive groups. The conceptual interpretation on the occurrence of classical Deligne-Lusztig varieties was given; however, a large part of the work in \cite{GH} was still obtained by brute force. 

\subsection{Some equivalent conditions} From the Lie-theoretic point of view, one would like to consider not only the maximal parahoric subgroups, but all parahoric subgroups; and one would like to have a conceptual understanding on the following question: 

When and why is $X(\mu, b_0)_{K'}$ naturally a union of classical Deligne-Lusztig varieties? 

This was finally achieved in joint work with G\"ortz and Nie \cite{GHN} as follows

\begin{theorem}\label{GHN-nice}
Assume that $G$ is simple, $\mu$ is a dominant coweight of $G$ and $K'$ is a parahoric subgroup. Then the following conditions are equivalent:

\begin{itemize}
\item For basic $[b_0] \in B(G, \mu)$, $X(\mu, b_0)_{K'}$ is naturally a union of classical Deligne-Lusztig varieties;

\item For any nonbasic $[b] \in B(G, \mu)$, $\dim X(\mu, b)_{K'}=0$;

\item The pair $(\mu, b)$ is Hodge-Newton decomposable for any nonbasic $[b] \in B(G, \mu)$; 

\item The coweight $\mu$ is minute for $G$. 
\end{itemize}
\end{theorem}

Here the minute condition is an explicit combinatorial condition on the coweight $\mu$. For quasi-split groups, it means that for any $\s$-orbit $\mathcal O$ on the set of simple roots, we have $\sum_{i \in \mathcal O} \<\mu, \omega_i\> \le 1$. For non quasi-split groups, the condition is more involved and we refer to \cite[Definition 2.2]{GHN2} for the precise definition. It is also worth mentioning that it is not very difficult to classify the pairs $(G, \mu)$ with the minute condition. In \cite[Theorem 2.5]{GHN2}, a complete list of the cases is obtained, where $X(\mu, b_0)_{K'}$ is naturally a union of classical Deligne-Lusztig varieties.

Fargues and Rapoport conjectured that for $p$-adic period domains, the weakly admissible locus coincides with the admissible locus if and only if the pair $(\mu, b)$ is Hodge-Newton decomposable for any nonbasic $[b] \in B(G, \mu)$ (cf. \cite[Conjecture 0.1]{GHN2}). This conjectured is established in a very recent preprint \cite{CFS} by Chen, Fargues and Shen. 

\subsection{Further remarks} From the Lie-theoretic point of view, there are some quite striking new features in Theorem \ref{GHN-nice}:

\begin{enumerate}
\item The relations between the variety $X(\mu, b)_{K'}$ for the basic $\s$-conjugacy class and for nonbasic $\s$-conjugacy classes;

\item The relation between the condition that $X(\mu, b_0)_{K'}$ has a simple description and the Hodge-Newton decomposability condition;

\item The existence of a simple description of $X(\mu, b_0)_{K'}$ is independent of the parahoric subgroup $K'$. 
\end{enumerate}

Note that part (1) and part (2) are new even for the specific Shimura varieties with hyperspecial level structure considered in the previous works. Part (3) is the most mysterious one. In \cite{GHN2}, we state that ``We do not see any reason why this independence of the parahoric could be expected a priori, but it is an interesting parallel with the question when the weakly admissible and admissible loci in the rigid analytic period domain coincide.''

For applications to number theory, one needs to consider the fraction field of the Witt ring instead of the formal Laurent series field $\overline \BF_q((\e))$. In that setting, we have a similar, but weaker result, namely,  $X(\mu, b_0)_{K'}$ is naturally a union of classical Deligne-Lusztig varieties as perfect schemes. It is expected that the structural results hold without perfection, as indicated in the special cases established in the papers mentioned in \S \ref{8.1}. 

\section{Some applications to Shimura varieties} In the last subsection, we give a very brief discussion of some applications to arithmetic geometry.

\subsection{Some characteristic subsets} The study of some characteristic subsets in the special fiber of a Shimura variety is a central topic in arithmetic geometry. We mention the \emph{Newton strata},  the \emph{Ekedahl-Oort strata} for the hyperspecial level structure and the \emph{Kottwitz-Rapoport strata} for the Iwahori level structure. Concerning these stratifications, there are many interesting questions one may ask, e.g. which strata are nonempty, what is the relation between these various stratifications, etc.. These questions have been intensively studied in recent years and there is a large body of literature on these questions. Among them, we mention the work of Viehmann and Wedhorn \cite{V-W} on the nonemptiness of Newton strata and Ekedahl-Oort strata for PEL type Shimura varieties with hyperspecial level structure, the work of Kisin, Madapusi and Shin \cite{KMS} on the nonemptiness of the basic Newton stratum, the work of Hamacher \cite{Ham} on the closure relation between Newton strata, and the work of Wedhorn \cite{We} and Moonen \cite{Mo} on the density of the $\mu$-ordinary locus (i.e. the Newton stratum corresponding to $\e^\mu$). We refer to \cite[Introduction]{HR} and \cite{Vx} for more references.

\subsection{An axiomatic approach} In the works mentioned above, both algebro-geometric and Lie-theoretic methods are involved, and are often mixed together. 

In joint work with Rapoport \cite{HR}, we purposed an axiomatic approach to the study of these characteristic subsets in a general Shimura variety. We formulated five axioms, based on the existence of integral models of Shimura varieties (which have been established in various cases by the work of Rapoport and Zink \cite{RZ}, Kisin and Pappas \cite{KP}), the existence of the following commutative diagram and some compatibility conditions:
\[
\xymatrix{& & {K'} \backslash G/{K'} \\  Sh_{K'}  \ar[r]^-{\Upsilon_{K'}} \ar@/^1pc/[urr]^{ \lambda_{K'} } \ar@/_1pc/[drr]_{ \delta_K } & G/{K'}_\sigma \ar[ur]_{\ell_K} \ar[dr]^{d_{K'}} & \\ & & B(G)}.
\]
Here $K'$ is a parahoric subgroup, $Sh_{K'}$ is the special fiber of a Shimura variety with $K'$ level structure, and $G/K'_\sigma$ is the set-theoretic quotient of $G$ by the $\s$-conjugation action of $K'$. 

As explained in \cite[\S 6.2]{GHN2}, affine Deligne-Lusztig varieties are involved in the diagram in an essential way, via the bijection \[
J_b\backslash X(\mu, b)_{K'} \isoarrow
d_{K'}^{-1}([b]) \cap \ell_{K'}^{-1}(K' \Adm(\{\mu\}) K').
\] 

\subsection{Some applications and current status of the axioms} 

It is shown in \cite{HR} that under those axioms, the Newton strata, the Ekedahl-Oort strata, the Kottwitz-Rapoport strata, and the Ekedahl-Kottwitz-Oort-Rapoport strata discussed in \S \ref{par}, are all nonempty in their natural range. Furthermore, under those axioms several relations between these various stratifications are also established in \cite{HR}. 

Following \cite{HR}, Shen and Zhang in \cite{SZ} studied the geometry of good reductions of Shimura varieties of abelian type. They established basic properties of these characteristic subsets, including nonemptiness, closure relations and dimension formula and some relations between these stratifications.

In joint work with Nie \cite{HN3}, based on the framework of \cite{HR}, we studied the density problem of the $\mu$-ordinary locus. Under the axioms of \cite{HR} we gave several explicit criteria on the density of the $\mu$-ordinary locus. 

Algebraic geometry is essential in the verification of these axioms. For PEL type Shimura varieties associated to unramified groups of type A and C and to odd ramified unitary groups, the axioms are verified in joint work with Zhou \cite{HZ}. For Shimura varieties of Hodge type, most of the axioms are verified in recent work of Zhou \cite{Zh}. 





\begin{thebibliography}{99}

\bibitem{BS}
B. Bhatt and P. Scholze, \emph{Projectivity of the Witt vector Grassmannian}, Invent. math. \textbf{209} (2017), 329--423.

\bibitem{BT}
F.~Bruhat and J.~Tits, \emph{Groupes r\`eductifs sur un corps local, I}, Publ. Math. IHES 41 (1972), 5--276.

\bibitem{Ch}
C.-L. Chai, \emph{Newton polygons as lattice points}, Amer. J. Math. 122 (5) (2000), 967--990.

\bibitem{CFS}
M. Chen, L. Fargues and X. Shen, \emph{On the structure of some $p$-adic period domains}, arXiv:1710.06935. 

\bibitem{CKV}
M. Chen, M. Kisin and E. Viehmann, \emph{Connected components of affine Deligne-Lusztig varieties in mixed
characteristic}, Compositio Math. \textbf{151} (2015), 1697--1762.

\bibitem{CN}
L. Chen and S. Nie, \emph{Connected components of closed affine Deligne-Lusztig varieties}, arXiv:1703.02476.

\bibitem{CH}
D. Ciubotaru and X. He, \emph{Cocenters and representations of affine Hecke algebra}, J. Eur. Math. Soc. \textbf{19} (2017), 3143--3177. 

\bibitem{DL}
P. Deligne and G. Lusztig, {\em Representations of reductive groups over finite fields}, Ann. of Math. (2) \textbf{103} (1976), 103--161.

\bibitem{DTX}
D. Helm, Y. Tian, L. Xiao, \emph{Tate cycles on some unitary Shimura varieties mod $p$}, Algebra Number Theory \textbf{11} (2017), 2213--2288.

\bibitem{Gai}
D. Gaitsgory, \emph{Construction of central elements in the affine Hecke algebra via nearby cycles}, Invent. Math. \textbf{144} (2001), 253--280.

\bibitem{Ga}
Q. Gashi, On a conjecture of Kottwitz and Rapoport, Ann. Sci. \'Ecole Norm. Sup. (4) \textbf{43} (2010), no. 6, 1017--1038. 

\bibitem{GP}
M.~ Geck and G.~ Pfeiffer, \emph{On the irreducible characters of {H}ecke algebras}, Adv. Math. \textbf{102} (1993), 79--94.

\bibitem{GHKR1} 
U.~G\"{o}rtz, T.~Haines, R.~Kottwitz, D.~Reuman, {\em Dimensions of some affine Deligne-Lusztig varieties}, Ann. sci. \'Ecole Norm. Sup. \textbf{39} (2006), 467--511.

\bibitem{GHKR2} 
U.~G\"{o}rtz, T.~Haines, R.~Kottwitz, D.~Reuman, {\em Affine Deligne-Lusztig varieties in affine flag varieties}, Compos. Math.\textbf{146} (2010), 1339-1382.

\bibitem{GH}
U. G\"ortz and X. He, \emph{Basic loci in Shimura varieties of Coxeter type}, Cambridge J. Math.~3 (2015), no.~3, 323--353.

\bibitem{GHN}
U. G\"ortz, X. He and S. Nie, \emph{$P$-alcoves and nonemptiness of affine Deligne-Lusztig varieties}, Ann. Sci. \`Ecole Norm. Sup. 48 (2015), 647--665.

\bibitem{GHN2}
U. G\"ortz, X. He and S. Nie, \emph{Fully Hodge-Newton decomposable Shimura varieties}, arXiv: 1610.05381.

\bibitem{GY}
U. G\"ortz and C.-F. Yu, \emph{Supersingular Kottwitz-Rapoport strata and Deligne-Lusztig varieties}, J. Inst. Math. Jussieu \textbf{9} (2010), 357--390. 

\bibitem{HH}
T. Haines and X. He, \emph{Vertexwise criteria for admissibility of alcoves}, Amer. J. Math. \textbf{139} (2017), 769--784.

\bibitem{HN} T.~Haines, B.~C.~Ng\^{o}, {\em Alcoves associated to special fibers of local models}, Amer. J. Math. {\bf 124} (2002), 1125-1152.

\bibitem{Ha1}
P.~Hamacher, \emph{The dimension of affine Deligne-Lusztig varieties in the affine Grassmannian of unramified groups}, Int. Math. Res. Notices \textbf{23} (2015), 12804--12839. 

\bibitem{Ham}
P.~Hamacher, \emph{The geometry of Newton strata in the reduction modulo $p$ of Shimura varieties of PEL type},  Duke Math. J. 164 (2015), 2809--2895. 

\bibitem{He07}
X.~He, \emph{The G-stable pieces of the wonderful compactification}, Trans. Amer. Math. Soc., 359 (2007), 3005--3024.

\bibitem{He072}
X.~He,  \emph{Minimal length elements in some double cosets of {C}oxeter groups}, Adv. Math. 215 (2007), no.~2, 469--503.

\bibitem{He11}
X. He, {\em Closure of Steinberg fibers and affine Deligne-Lusztig varieties}, Int. Math. Res. Notices 14 (2011), 3237--3260.

\bibitem{He14}
X.~He, \emph{Geometric and homological properties of affine Deligne-Lusztig varieties}, Ann. Math. \textbf{179} (2014), 367--404.  

\bibitem{He16}
X. He, \emph{Kottwitz-Rapoport conjecture on unions of affine Deligne-Lusztig varieties}, Ann. Sci. \`Ecole Norm. Sup. 49 (2016), 1125--1141.

\bibitem{CDM}
X. He, \emph{Hecke algebras and $p$-adic groups}, Current developments in mathematics 2015, 73--135, Int. Press, Somerville, MA, 2016.

\bibitem{HN14}
X.~He and S.~ Nie, \emph{Minimal length elements of extended affine Weyl group}, Compos. Math. 150 (2014), no. 11, 1903--1927.

\bibitem{HNx}
X. He and S. Nie, \emph{On the acceptable elements}, IMRN \textbf{2018} (2018), 907--931.

\bibitem{HN3}
X. He and S. Nie, \emph{On the $\mu$-ordinary locus of a Shimura variety}, Adv. Math. \textbf{321} (2017), 513--528.

\bibitem{HR}
X. He and M. Rapoport, \emph{Stratifications in the reduction of Shimura varieties}, Manuscripta Math. 152 (2017), 317--343.

\bibitem{HZ}
X. He and R. Zhou, \emph{On the connected components of affine Deligne-Lusztig varieties}, arXiv:1610.06879.

\bibitem{HP2}
B.~Howard and G.~Pappas, \emph{On the supersingular locus of the $GU(2,2)$ Shimura variety}, Algebra Number Theory 8 (2014), no. 7, 1659--1699. 

\bibitem{HP}
B.~Howard, G.~Pappas, \emph{Rapoport-Zink spaces for spinor groups}, Compos. Math. \textbf{153} (2017), 1050--1118.

\bibitem{IM}
N.~Iwahori and H.~Matsumoto, \emph{On some Bruhat decomposition and the structure of the Hecke rings of $p$-adic Chevalley groups}, Publ. Math. IHES 25 (1965), 5--48.

\bibitem{Ka}
H. Katz, \emph{Slope filtration of $F$-crystals}, Ast\'erisque 63 (1979), 113--163.

\bibitem{Ki}
M. Kisin, \emph{Mod-$p$ points on Shimura varieties of abelian type}, J. Amer. Math. Soc. 30 (2017), 819--914.

\bibitem{KMS}
M.~Kisin, K. Madapusi, S.-W. Shin, \emph{Honda-Tate theory for Shimura varieties of Hodge type}, in preparation.

\bibitem{KP}
M. Kisin and G. Pappas, \emph{Integral models of Shimura varieties with parahoric level structure}, arXiv:1512.01149, to appear in Publ. Math. IHES.

\bibitem{K1} 
R.~Kottwitz, \emph{Isocrystals with additional structure}, Compositio Math.  \textbf{56} (1985), 201--220.

\bibitem{K2}
R.~Kottwitz, \emph{Isocrystals with additional structure. {II}}, Compositio Math. \textbf{109} (1997), 255--339.

\bibitem{K3}
R.~Kottwitz, \emph{On the Hodge-Newton decomposition for split groups}, IMRN 26 (2003), 1433--1447.

\bibitem{K4}
R.~Kottwitz, \emph{Dimensions of Newton strata in the adjoint quotient of reductive groups}, Pure Appl. Math. Q. \textbf{2} (2006), 817--836.

\bibitem{KR}
R. Kottwitz, M. Rapoport, {\em Minuscule alcoves for $GL_n$ and $GSP_{2n}$}, Manuscripta Math. {\bf 102}, no.4, (2000), 403-428.

\bibitem{KR03}
R.~ Kottwitz and M.~ Rapoport, \emph{On the existence of F -crystals}, Comment. Math. Helv., \textbf{78} (2003), 153--184.

\bibitem{KR2} 
S. Kudla and M. Rapoport, \emph{Special cycles on unitary Shimura varieties I. Unramified local theory}, Invent. Math. 184 (3) (2011), 629--682.

\bibitem{LR}
R. Langlands, M. Rapoport, Shimuravariet\"aten und Gerben, J. Reine Angew. Math 378 (1987), 113--220.

\bibitem{LZ}
C. Li and Y. Zhu, \emph{Remarks on the arithmetic fundamental lemma}, Algebra Number Theory \textbf{11} (2017), 2425--2445.

\bibitem{L-orange}
G.~Lusztig, \emph{Characters of reductive groups over a finite field}, Annals of Mathematics Studies, 107. Princeton University Press, Princeton, NJ, 1984.

\bibitem{L1}
G.~Lusztig, {\em Parabolic character sheaves}, I, II, Mosc. Math. J. 4 (2004), no.~1, 153--179; no.~4, 869--896.

\bibitem{L1-3}
G. Lusztig, \emph{From conjugacy classes in the Weyl group to unipotent classes}, I--III, Represent.Th. 15(2011), 494-530; 16(2012), 189--211, 450-488.

\bibitem{Lu-Bruhat}
G. Lusztig, \emph{Bruhat decomposition and applications}, arXiv:1006.5004.

\bibitem{L99}
G. Lusztig, \emph{Parabolic character sheaves}, III, Mosc. Math. J. 10 (2010), 603--609.

\bibitem{Ma}
B.~Mazur, \emph{Frobenius and the Hodge filtration (estimates)}, Ann. of Math. (2) 98 (1973), 58--95.

\bibitem{MST}
E. Mili\'cevi\'c, P. Schwer and A. Thomas, \emph{Dimensions of affine Deligne-Lusztig varieties: a new approach via labeled folded alcove walks and root operators}, arXiv:1504.07076. 

\bibitem{Mo}
B. Moonen, \emph{Serre-Tate theory for moduli spaces of PEL-type}, Ann. Sci. \`Ecole Norm. Sup. 37 (2004), 223--269.

\bibitem{Nie}
S.~Nie, \emph{Connected components of closed affine Deligne-Lusztig varieties in affine Grassmannians}, arXiv:1511.04677.

\bibitem{PR}
G. Pappas and M. Rapoport, \emph{Twisted loop groups and their affine flag varieties}, Adv. Math. \textbf{219} (2008), 118--198, with an appendix by T. Haines and Rapoport.

\bibitem{PWZ}
R. Pink, T. Wedhorn, and P. Ziegler, \emph{Algebraic zip data} Doc. Math. \textbf{16} (2011), 253--300.

\bibitem{RH}
E. A. Rad and U. Hartl, \emph{Langlands-Rapoport Conjecture Over Function Fields}, arXiv:1605.01575.

\bibitem{Ra}
M.~Rapoport, \emph{A guide to the reduction modulo {$p$} of {S}himura varieties}, Ast\'erisque (2005), no.~298, 271--318.

\bibitem{RR}
M.~Rapoport and M.~Richartz, \emph{On the classification and specialization of $F$-isocrystals with additional structure}, Compositio Math. 103 (1996), no. 2, 153--181. 

\bibitem{RTW} M.~Rapoport, U.~Terstiege, S.~Wilson, \emph{The supersingular locus of the Shimura variety for $GU(1,n-1)$ over a ramified prime}, Math. Z. 276 (2014), 1165--1188.

\bibitem{RTZ} M.~Rapoport, U.~Terstiege, W.~Zhang, \emph{On the Arithmetic Fundamental Lemma in the minuscule case}, Compos.~Math. \textbf{149} (2013), 1631--1666.

\bibitem{RZ} 
M.~Rapoport, Th.~Zink, \emph{Period spaces for $p$-divisible groups}, Annals of Mathematics Studies, 141. Princeton University Press, Princeton, NJ, 1996.

\bibitem{SZ}
X. Shen and C. Zhang, \emph{Stratifications in good reductions of Shimura varieties of abelian type}, arXiv:1707.00439.

\bibitem{TX} Y.~Tian, L.~Xiao, \emph{On Goren-Oort stratification for quaternionic Shimura varieties}, Compos. Math. 152 (2016), 2134--2220.

\bibitem{TX2} Y.~Tian, L.~Xiao, \emph{Tate cycles on quaternionic Shimura varieties over finite fields}, arXiv:1410.2321 (2014).

\bibitem{Vi06} 
E.~Viehmann. {\em The dimension of affine Deligne-Lusztig varieties},  Ann. Sci.
\'Ecole Norm. Sup. (4) \textbf{39} (2006), 513--526.

\bibitem{Vi08}
E.~Viehmann, \emph{Connected components of affine Deligne-Lusztig varieties},
Math. Ann. \textbf{340} (2008), 315--333.

\bibitem{Vi2}
E.~Viehmann, \emph{Truncations of level $1$ of elements in the loop group of a reductive group}, Ann. of Math. (2) 179 (2014), 1009--1040. 

\bibitem{Vx}
E.~Viehmann, \emph{On the geometry of the Newton stratification}, arXiv:1511.03156.

\bibitem{V-W}
E. Viehmann, T. Wedhorn, \emph{Ekedahl-Oort and Newton strata for Shimura varieties of PEL type}, Math. Ann. 356 (2013), 1493--1550. 

\bibitem{VW} 
I.~Vollaard, T.~Wedhorn, \emph{The supersingular locus of the Shimura variety of GU(1,n-1) II}, Invent.~Math.~\textbf{184} (2011), 591--627.

\bibitem{We}
T. Wedhorn, \emph{Ordinariness in good reductions of Shimura varieties of PEL-type}, Ann. Sci. \`Ecole Norm. Sup. (4), 32 (1999), 575--618.

\bibitem{Wi} 
J.-P. Wintenberger, \emph{Existence de {$F$}-cristaux avec structures suppl\'ementaires}, Adv. Math. \textbf{190} (2005), 196--224.

\bibitem{Zhang}
W.~Zhang, \emph{On arithmetic fundamental lemmas}, Invent. Math., \textbf{188} (2012), 197--252.

\bibitem{Zh}
R. Zhou, \emph{Mod-$p$ isogeny classes on Shimura varieties with parahoric level structure}, arXiv: 1707.09685.

\bibitem{Z1}
X.~Zhu, \emph{On the coherence conjecture of Pappas and Rapoport}, Ann. of Math. (2) \textbf{180} (2014), 1--85.

\bibitem{Z2}
X.~Zhu, \emph{Affine Grassmannians and the geometric Satake in mixed characteristic}, Ann. of Math. (2) \textbf{185} (2017),  403--492.

\end{thebibliography}
\end{document}